\documentclass[journal]{IEEEtran}
\usepackage[pdftex]{graphicx}
\usepackage{amsfonts}
\usepackage{setspace}
\usepackage{wrapfig}
\usepackage{color}
\usepackage{cancel}
\usepackage{amsthm}
\usepackage{url}
\usepackage{slashbox}
\usepackage{enumerate}
\usepackage{cite}
\usepackage{amssymb}
\usepackage[cmex10]{amsmath}
\newcommand{\field}[1]{\mathbb{#1}}                                   
\newcommand{\prn}[1]{\left(#1\right)}                                 
\newcommand{\Prn}[1]{\left[#1\right]}                                 
\newcommand{\dff}[2]{\displaystyle\frac{d#1}{d#2}}                    
\newcommand{\dlims}[2]{\displaystyle\lim_{#1}{\ #2}}          				
\newcommand{\dint}[3]{\displaystyle\int_{#1}^{#2}{#3}}                
\newcommand{\dsum}[3]{\displaystyle\sum_{#1}^{#2}{#3}}                
\newcommand{\abs}[1]{\left|#1\right|}                                 
\newcommand{\eltwo}[1]{\left\|#1\right\|_2}                           
\newcommand{\norm}[1]{\left\|#1\right\|}                              
\newcommand{\sarray}[1]{\begin{array}{l}#1\end{array}}								
\newcommand{\amatrix}[1]{\begin{bmatrix}#1\end{bmatrix}}							
\newcommand{\sign}[1]{\ \mathrm{sign}{\left(#1\right)}}								
\newcommand{\set}[1]{\displaystyle\left\{#1\right\}}									
\newtheorem{theorem}{Theorem}
\newtheorem{lemma}{Lemma}

\newtheorem{definition}{Definition}

\theoremstyle{remark}
\newtheorem*{remark}{Remark}
\newcommand{\G}{\Gamma}
\newcommand{\Gc}{\Gamma^c}

\newcommand{\PP}{\Phi_{\Gamma}^T\Phi_{\Gamma}}
\newcommand{\PcP}{\Phi_{\Gamma^c}^T\Phi_{\Gamma}}

\begin{document}
%
\title{Convergence and Rate Analysis of\\Neural Networks for Sparse Approximation}
%
%
%

\author{Aur\`ele~Balavoine,~\IEEEmembership{Student Member,~IEEE,}
		Justin~Romberg,~\IEEEmembership{Senior member,~IEEE,}
        and~Christopher~J.~Rozell,~\IEEEmembership{Member,~IEEE,}
\thanks{The authors are with the School of Electrical and Computer Engineering, Georgia Institute of Technology, Atlanta,
GA, 30332-0250 USA.}\thanks{e-mail: \{aurele.balavoine,jrom,crozell\}@gatech.edu,}}

\maketitle

\begin{abstract}
We present an analysis of the Locally Competitive Algorithm (LCA), a Hopfield-style neural network that efficiently solves sparse approximation problems (e.g., approximating a vector from a dictionary using just a few non-zero coefficients).  This class of problems plays a significant role in both theories of neural coding and applications in signal processing. However, the LCA lacks analysis of its convergence properties and previous results on neural networks for nonsmooth optimization do not apply to the specifics of the LCA architecture.
We show that the LCA has desirable convergence properties, such as stability and global convergence to the optimum of the objective function when it is unique. Under some mild conditions, the support of the solution is also proven to be reached in finite time.  Furthermore,  some restrictions on the problem specifics allow us to characterize the convergence rate of the system by showing that the LCA converges exponentially fast with an analytically bounded convergence rate.  We support our analysis with several illustrative simulations.
\end{abstract}

\makeatletter{\renewcommand*{\@makefnmark}{}
\footnotetext{This work was partially supported by NSF grant CCF-0905346 and NIH grant R01-EY019965.  Preliminary versions of portions of this work were presented in~\cite{balavoine_global_2011}.\\
This paper was published in IEEE trans. Neural Netw., vol. 23, no. 9, pp. 1377--1389, Sept. 2012.\\
The statement and proof of Theorem 3 have been corrected in this version.}\makeatother}

\section{Introduction}
\label{sec:intro}

\IEEEPARstart{S}{parse} approximation has generated substantial interest in a wide range of research communities over the last two decades, including signal processing, machine learning, statistics and computational neuroscience (e.g., see~\cite{elad_role_2010,olshausen_sparse_2004} and references therein). Specifically, sparse approximation involves solving an optimization problem to represent a signal using just a few atoms from some specified (possibly overcomplete) dictionary.  In addition to describing compelling models of neural coding for sensory information~\cite{olshausen_emergence_1996,olshausen_sparse_1997}, this approach has led to state of the art results in many types of inverse problems.  One example of a regime that has leveraged this type of signal model is the theory of Compressed Sensing (CS)~\cite{candes_robust_2006,donoho_compressed_2006}.  In this domain, highly undersampled signals are recovered by solving a sparse approximation problem, thereby shifting the burden of data acquisition from the front-end sensor to a computationally intensive back end.

Because of the increasing interest in the sparse approximation paradigm, significant effort has been made to design efficient algorithms for solving (or approximately solving) least squares optimization problems that are regularized with a sparsity-inducing penalty (e.g.,~\cite{kim_interior-point_2007,figueiredo_gradient_2007,daubechies_iterative_2004}).  However, all of these algorithms are developed to run on digital computers in discrete time, which are both implausible for neural systems and suffer from several drawbacks as engineering approaches.  In particular, the computational time required by these algorithms presents a barrier to real-time signal processing applications with high dimensional signals at significant data rates.  Specifically, digital algorithms tend to have storage requirements and convergence times that scale unfavorably with the dimensionality of the signals being approximated.  Additionally, the power consumption of digital solutions can be prohibitive for many applications.  Given these considerations, a fast and low-power method for solving sparse approximation problems in a parallel architecture would be valuable both for engineering systems and viable models of neural coding.

Analog neural networks have long been proposed for solving optimization problems~\cite{cichocki_neural_1993}, with an early example being Hopfield's pioneering results~\cite{hopfield_neural_1982} using networks to solve linear optimization problems.   Analog neural networks offer several potential advantages over comparable digital algorithms, including their ability to be implemented in analog architectures that are highly parallel, fast and power efficient. Recent advances in VLSI reconfigurable analog chips~\cite{twigg_configurable_2009} make the design of such systems more feasible and affordable than has often been true in the past. 
%
%
A recent neural network architecture, called a Locally Competitive Algorithm (LCA)~\cite{rozell_sparse_2008}, has been proposed to solve the types of nonsmooth optimization problems that come up in sparse approximation.  These Hopfield-style networks appear to efficiently solve a whole family of sparse approximation problems by incorporating ideal thresholding nonlinearities in the network dynamics. The results in~\cite{rozell_sparse_2008} provide encouraging evidence that biological systems and engineering applications can use neural networks to solve these important sparse approximation problems.

However, despite encouraging evidence of its performance, the LCA lacks strong convergence guarantees and estimates of the convergence rate. Furthermore, the specifics of the LCA architecture violate many of the assumptions used in previous analysis, making the extensive literature on neural network convergence inapplicable for this system. For example, in previous work the activation function is often assumed to be linear~\cite{hespanha_uniform_2004}, piecewise linear~\cite{liu_simplified_2006}, or nonlinear but increasing and bounded~\cite{yang_exponential_1994,liang_recurrent_2000,forti_absolute_2004,lu_convergence_2008}. In contrast, the LCA activation function is nonlinear and unbounded. Other analyses assume an interconnection matrix which is positive definite~\cite{xia_recurrent_2004,ferreira_support_2006,liu_simplified_2006} or non-singular~\cite{forti_new_1995,li_dynamical_2009}, while the LCA interconnection matrix may have negative eigenvalues as well as a non-trivial nullspace (due to the approximation dictionary being overcomplete). Other analyses consider a nonsmooth objective function with constraints, but only shows convergence when the constraints are nonzero and convex~\cite{forti_generalized_2004,forti_convergence_2006}. 
Finally, other relevant analyses assume the objective function is convex~\cite{xue_subgradient-based_2008} or piecewise linear~\cite{liu_finite-time_2011}, while the work in~\cite{lin_stability_2009} focuses on controllability of the network path.

The main contributions of this paper are to present a formal analysis of the LCA network architecture and its convergence properties, despite 1) an activation function that is nonsmooth and not necessarily bounded or increasing, and 2) a potentially singular interconnection matrix.  Section~\ref{sec:conv} contains our first main result, which states that the fixed points of the LCA network correspond to critical points of the objective function. In the special case where the objective is convex, this set coincides with the global minima of the objective. In addition, we show that the network is globally stable and that the outputs are quasi-convergent, in the sense that they get infinitely close to a set of fixed points. Finally, in the case where the objective function has isolated minima, we show that the LCA converges to the solution of the sparse approximation problem for any initial point (i.e., a much stronger condition than just the non-increasing property of the energy function that was shown in previous work~\cite{rozell_sparse_2008}).  This section also shows that the LCA is well-behaved in that it converges in a finite number of switches (i.e., nodes crossing above or below threshold).  Section~\ref{sec:rate} expands on these general results to show that, under additional mild conditions on the problem specifics, the LCA actually converges \emph{exponentially fast} to the solution.  Furthermore, we give an analytic expression for this convergence rate that depends on the properties of the detailed approximation problem.  Finally, Section~\ref{sec:sim} presents simulation results showing the correspondence of our analytic results with empirical observations of the network behavior.

\section{Background}
\label{sec:lca}

Before presenting our main results, in this section we briefly give a precise statement of the sparse approximation problems of interest, a description of the LCA architecture, and some preliminary observations on the LCA network dynamics that will be useful in the subsequent analysis.

\subsection{Sparse Approximation}
\label{ssec:approx}

As mentioned above, sparse approximation is an optimization program that seeks to find the approximation coefficients of a signal on a prescribed dictionary, using as few non-zero elements as possible. To fix notation, we denote the input signal by $y \in\mathbb{R}^M$, the unit-norm dictionary elements by the columns of the $M\times N$ matrix $\Phi =\left[\Phi_1, \ldots, \Phi_N\right]$, and the coefficients by $a \in\mathbb{R}^N$. Generally, $M\ll N$ (i.e., the approximation dictionary is overcomplete), and the problem of recovering $a$ from $y$ is underdetermined. While similar in spirit to the well-known winner-take-all problem~\cite{maass_computational_2000}, sparse approximation problems are generally formulated as the solution to an optimization program because this approach can often yield strong performance guarantees in specific applications (e.g., recovery in a CS problem).  In the most generic form, the objective function is the sum of a quadratic data fidelity term (i.e., mean squared error) and a regularization term that uses a sparsity-inducing cost function $C(\cdot)$:
\begin{equation}
\min\limits_{a}{\ V\left(a(t)\right)=\dfrac{1}{2}\left\|y-\Phi a\right\|^2_2 + \lambda \dsum{n=1}{N}{C(a_n)}},
\label{eq:lp}
\end{equation}
where the parameter $\lambda$ is a tradeoff between the two terms in the objective function.  The ``ideal'' sparse approximation problem has a cost function that simply counts the number of non-zero elements, resulting in a non-convex objective function that has many local minima~\cite{natarajan_sparse_1995}.

One of the most widely used programs from this family is known as Basis Pursuit De-Noising (BPDN) \cite{chen_atomic_2001}, which is given by the objective function:
\begin{equation}
\min\limits_{a}{\ \dfrac{1}{2}\left\|y-\Phi a\right\|^2_2 + \lambda \left\|a\right\|_1}.
\label{eq:l1}
\end{equation}
In BPDN, the $\ell_1$-norm $\left(\left\|a\right\|_1=\sum_n |a_n|\right)$ is used as a convex surrogate for the idealized counting norm. This program has gained in popularity as researchers have shown that, in many cases of interest, substituting the $\ell_1$ norm yields the same solution as using an idealized (and generally intractable) counting norm~\cite{donoho_sparse_2005}.  However, BPDN illustrates the canonical challenge of sparse approximation problems.  Despite being convex, the BPDN objective function contains a nonsmooth nonlinearity that makes it considerably more difficult than a classic least-squares problem.

In the context of computational neuroscience, sparse approximation has been proposed as a neural coding scheme for sensory information.  In one interpretation, programs such as BPDN can be viewed as Bayesian inference in a linear generative model with Gaussian noise and a prior with high kurtosis to encourage sparsity (e.g., the Laplacian prior in the case of BPDN)~\cite{olshausen_sparse_1997}.  Given the prevalence of probabilistic inference as a successful description of human perception~\cite{knill_bayesian_2004} and the theoretical benefits of sparse representations~\cite{olshausen_sparse_2004}, it has long been conjectured that sensory systems may encode stimuli via sparse approximation.  In fact, in classic results, sparse approximation applied to the statistics of natural stimuli in an unsupervised learning experiment has been shown to be sufficient to qualitatively and quantitatively explain the receptive field properties of simple cells in the primary visual cortex~\cite{olshausen_emergence_1996,rehn_network_2007} as well as the auditory nerve fibers~\cite{smith_efficient_2006}.  Only recently have there been proposals of efficient neural networks that could efficiently solve the necessary optimization problems to implement this type of encoding \cite{rozell_sparse_2008,rehn_network_2007,fischer_sparse_2007,perrinet_feature_2004}.

\begin{figure}
\centering
\includegraphics[width=3.2in, trim=0.5in 1.5in 0.5in 1.5in]{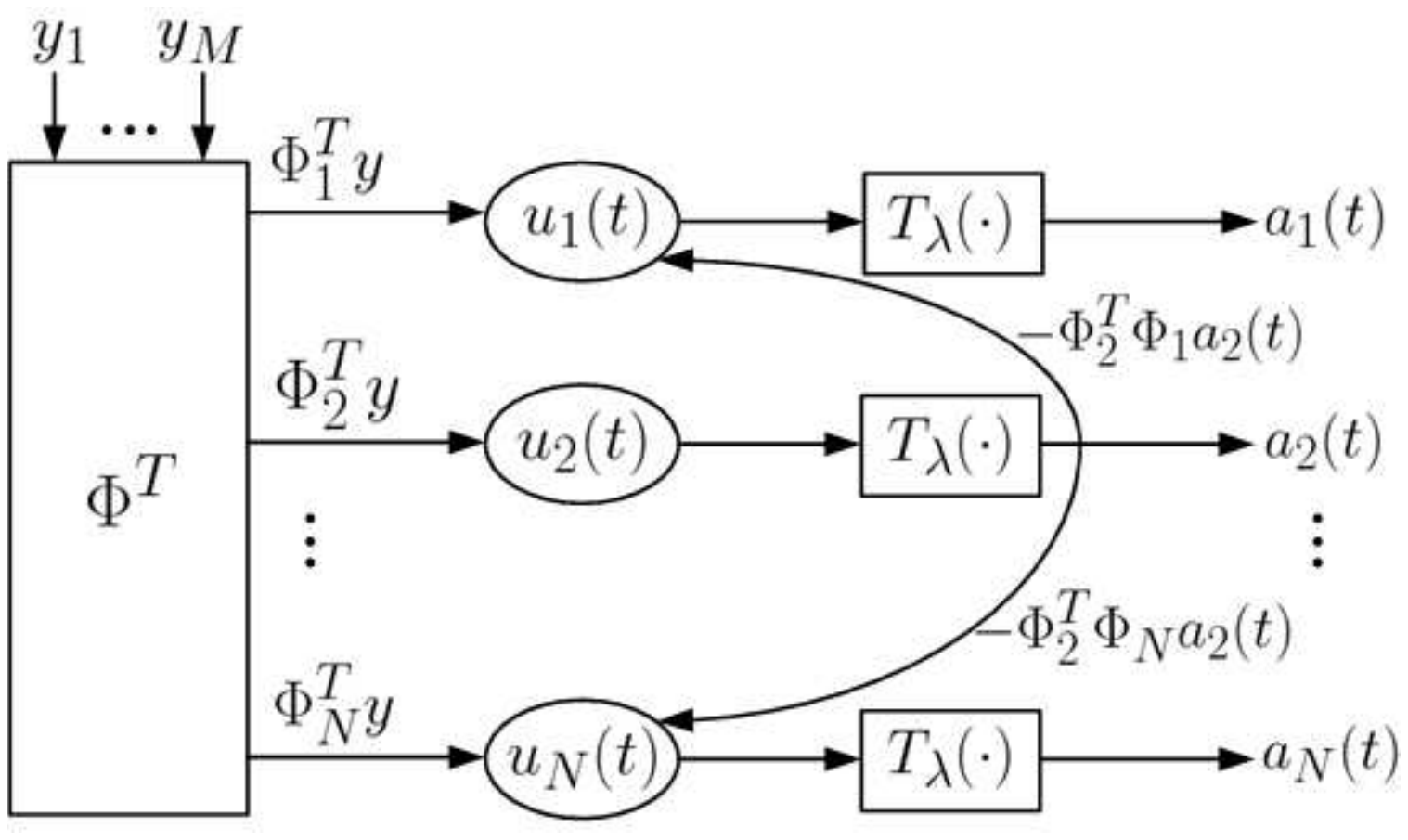}
\caption{Block diagram of the Locally Competitive Algorithm (LCA)~\cite{rozell_sparse_2008}, a Hopfield-style network architecture for solving sparse approximation problems.}
\label{fig:lca}
\end{figure}

\subsection{LCA structure}
\label{ssec:lca}

Our primary interest will be the LCA~\cite{rozell_sparse_2008}, an analog, continuous-time dynamical system that is a type of Hopfield-style network.  In particular, each node in the LCA network is characterized by the evolution of a set of internal state variables, $u_n(t)$ for $n=1, \ldots, N$, and uses a nonlinear activation function $T_{\lambda}(\cdot)$ to produce output variables $a_n(t)$ for $n=1, \ldots, N$. The activation function is typically a nonsmooth nonlinear function (such as a thresholding function) to induce sparsity in the outputs.  The dynamics of the internal state variables are governed by a set of coupled, nonlinear ordinary differential equations (ODEs):
\begin{equation} 
\begin{split}
	\tau \dot{u}(t)&=-u(t) - (\Phi^T \Phi - I) \ a(t) + \Phi^Ty\\
	a(t) &= T_{\lambda}(u(t))\\
\end{split} .
\label{eq:dyn}
\end{equation}
In this network, each node is associated with a single dictionary element $\Phi_n$, $n=1,\ldots,N$, and the node outputs will be shown to be the solution to the optimization problem of interest.  The architecture of a typical LCA is shown in \figurename~\ref{fig:lca}.

The inputs to the LCA network are feedforward connections computing the vector of driving synaptic inputs $\Phi^T y$, reflecting how well the signal $y$ matches each dictionary element.  The network also has recurrent inhibitory or excitatory connections between the nodes, modulated by weights corresponding to the interconnection matrix $G=\Phi^T\Phi-I$ (i.e., a modified Grammian matrix for the dictionary).  The more overlap there is between a pair of nodes (characterized by the inner product between their dictionary elements), the stronger the potential inhibition between those nodes.  While the modulating weight is symmetric between any pair of nodes, the total inhibition is not because it is also modulated by the activity of each individual node. Moreover, the matrix $G$ potentially has both negative and positive eigenvalues, as well as a non-trivial nullspace. This inhibition structure ensures that nodes that carry the same information will not become active at the same time, thus meeting the goals of sparse approximation. The time constant $\tau$ is dependent on the physical properties of the solver implementing the ODEs. For our purposes, we will often assume $\tau=1$ without loss of generality.

\label{ssec:actFunc}

\begin{figure}
\centering
\includegraphics[width=2.5in]{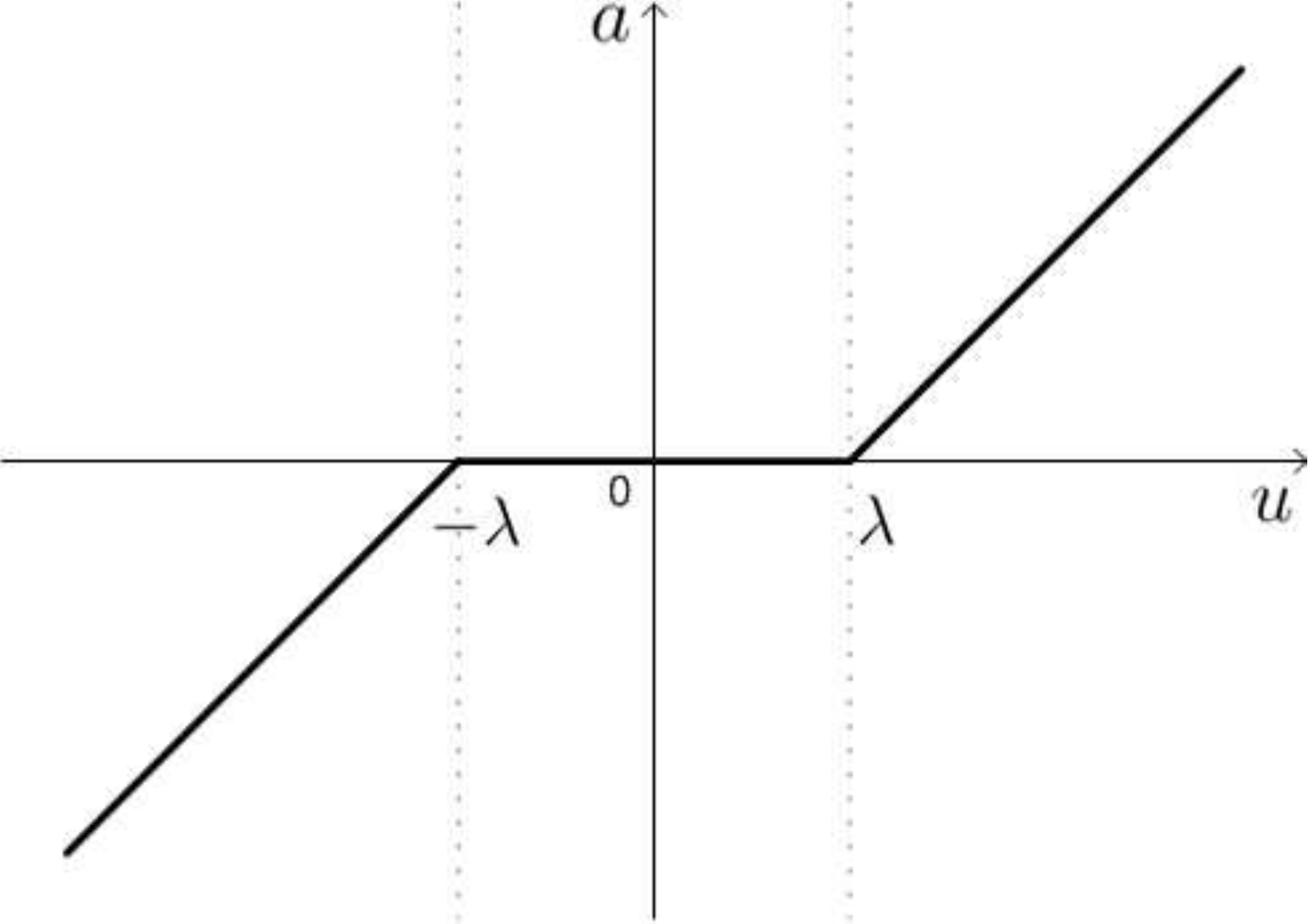}
\caption[Soft-Threshold]{The soft-thresholding activation function.  When this is used for the thresholding function $T_{\lambda}(\cdot)$, the LCA solves the popular BPDN optimization problem used in many sparse approximation applications.}
\label{fig:thresh}
\end{figure}

\begin{figure}
\centering
\includegraphics[width=2.5in, trim= 1in 0.5in 1in 0.5in]{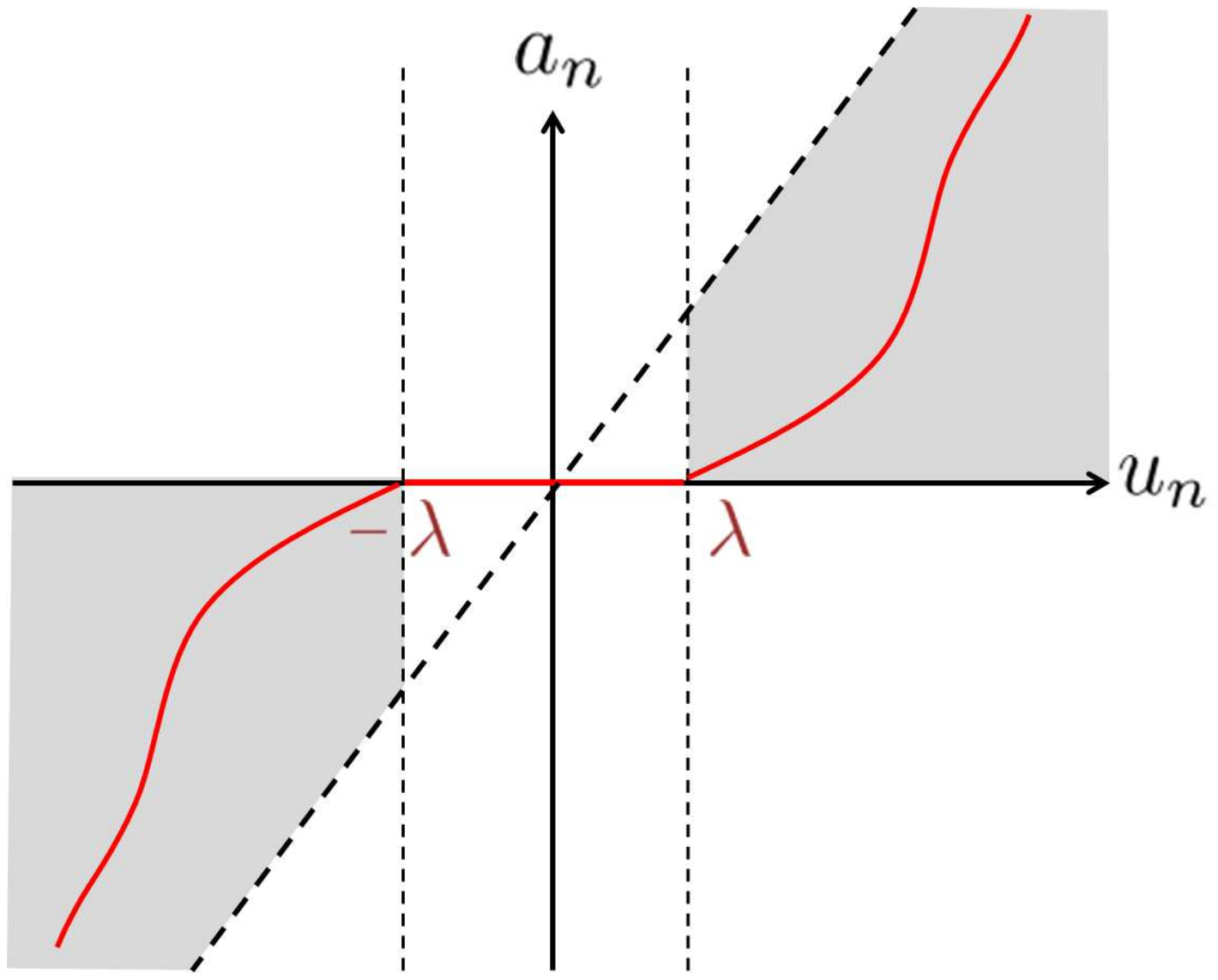}
\caption[Soft-Threshold]{Generic activation function satisfying conditions \eqref{eq:cond}. The area in gray represents where the activation function must lie in to satisfy conditions \eqref{eq:cond2} and \eqref{eq:cond3}.}
\label{fig:actFunc}
\end{figure}

It was shown in \cite{rozell_sparse_2008} that the objective function in~\eqref{eq:lp} is non-increasing along the LCA trajectory with the following relationship ($\forall a_n\in\field{R}$ such that $a_n\neq0$) between the cost penalty term $C(\cdot)$ and the activation function $T_{\lambda}(\cdot)$: 
\begin{equation}
\lambda\dff{C(a_n)}{a_n}=u_n-a_n=u_n-T_{\lambda}(u_n) .
\label{eq:corres}
\end{equation}
In the case of BPDN, the cost penalty is simply $C(x)=\abs{x}$ and the activation function obtained from \eqref{eq:corres} is the soft-thresholding function (\figurename~\ref{fig:thresh}) defined by:
$$a_n(t) = T_{\lambda}(u_n(t)) =\left\{
\begin{array}{ll}  0,& \abs{u_n(t)} \leq \lambda \\
u_n(t)-\lambda \sign{u_n(t)}, & \abs{u_n(t)} > \lambda 
\end{array}. \right. $$
This function often arises in connection to algorithms for minimizing the absolute value of the coefficients (e.g., see~\cite{bioucas-dias_new_2007}).
Generalizing from the soft-threshold, we focus on ``thresholding'' activation functions $T_{\lambda}(\cdot)$ of the form:
\begin{equation}
a_n(t) = T_{\lambda}(u_n(t)) =\left\{
\begin{array}{ll}  0,& \abs{u_n(t)} \leq \lambda \\
f(u_n(t)), & \abs{u_n(t)} > \lambda
\end{array}, \right.
\label{eq:thresh}
\end{equation}
where the function $f(\cdot)$ is a real-valued function defined and continuous on the domain $\mathcal{D}=\left(-\infty,-\lambda\right]\cup\left[\lambda,+\infty\right)$, differentiable on the interior of $\mathcal{D}$, and satisfying the conditions:
\begin{subequations}
\label{eq:cond}
        \begin{align}
                & f(-u_n)=-f(u_n), & & \forall u_n\in\mathcal{D} \mathrm{\ and\ } f(\lambda)=0 \label{eq:cond1} \\
                & f'(u_n)>0,       & & \forall u_n\in\mathcal{D} \label{eq:cond2} \\
                & f(u_n)\leq u_n,  & & \forall u_n\in\mathcal{D} \mathrm{\ s.t.\ } u_n\geq0 \label{eq:cond3} .
        \end{align}
\end{subequations} 
\begin{remark} Eq. \eqref{eq:cond1} ensures that $T_{\lambda}(\cdot)$ is continuous on all $\field{R}$ and consequently, that $u(t)$ and $a(t)$ are continuous with respect to time. This ensures that the cost function $C(\cdot)$ is differentiable everywhere except at the origin. Eq. \eqref{eq:cond2} makes $f(\cdot)$ a bijection on $\mathcal{D}$. Finally, Eq. \eqref{eq:cond1} and \eqref{eq:cond3} ensure, by a simple computation, that $C(\cdot)$ is positive and non-decreasing with the absolute value of the coefficient. Notice that the $\ell_1$-minimization objective function satisfies the three criteria in \eqref{eq:cond}, in addition to being convex (which is not required, but guarantees that the system will find a global minima).  In \figurename~\ref{fig:actFunc}, we plot a more generic stylized activation function that satisfies the conditions \eqref{eq:cond}.
\end{remark}

As shown in \figurename~\ref{fig:actFunc}, the activation function \eqref{eq:thresh} is composed of two operating regions. The first region corresponds to the case where the internal state $u_n$ is below the threshold $\lambda$, in which case the output $a_n$ is zero. We call the nodes in this region \emph{inactive nodes}. The second region corresponds to the case where the internal state $u_n$ is above threshold, in which case \eqref{eq:cond2} guarantees that the output $a_n$ is strictly increasing with $u_n$. We call the nodes that are above threshold \emph{active nodes}, and we denote by $\G$ the \emph{active set} (i.e., the set of indices corresponding to those nodes, $\G=\{n\in\Prn{1,N}: \abs{u_n(t)}>\lambda\}$). On the contrary, we denote the set of indices corresponding to nodes that are below threshold by $\Gc$ and call them the \emph{inactive set}. While these two sets do change with time as the network evolves, for the sake of readability we omit the dependence on time in the notation.

Due to the nonsmooth nature of the objective \eqref{eq:lp}, we require the generalized notion of a subgradient (see \cite{clarke_optimization_1987} for more details). Note that $V\left(a\right)$ is differentiable everywhere except at points $a$ that contain zeros, due to the non-differentiability of $C(\cdot)$ at the origin noted above for the cases of interest in sparse approximation. The subgradient $\partial V(\cdot)$ extends the classic notion of a gradient $\nabla V(\cdot)$ at those points of discontinuity. The properties in \eqref{eq:cond} allow us to use the following simple definition:
$$\partial V(a) = co\set{\lim\limits_{i\rightarrow\infty} \nabla V(a_i) : a_i\rightarrow a, a_i\notin\mathcal{S}, a_i\notin\Omega_V},$$
where $co$ is the convex hull, $\Omega_V$ is the set of point where $V$ fails to be differentiable and $\mathcal{S}$ is any set of Lebesgue measure $0$ in $\field{R}^N$. In other words, it is the smallest convex set containing the gradients of the function as we approach the discontinuity from any direction. Note that when $a$ is not a point of discontinuity, this simply reduces to $\nabla V(a)$. 

\subsection{LCA Node Dynamics}

Because nodes can cross threshold and go from the inactive set to the active set (and vice versa), the LCA can be thought of as a switched system \cite{decarlo_perspectives_2000} where the ODEs change at each switching time (i.e., when a node crosses above or below threshold).  In between two switching times, the active set $\G$ and inactive set $\Gc$ are fixed. Therefore, the dynamics on the two sets can be considered separately to facilitate analysis. In the following, $\Phi_{\mathcal{T}}$ is the matrix composed of the columns of $\Phi$ indexed by the set $\mathcal{T}$. Similarly, $u_{\mathcal{T}}$ and $a_{\mathcal{T}}$ refer to the elements in the original vectors indexed by ${\mathcal{T}}$. We also denote by $\{t_k\}_{\{k\in\field{N}\}}$ the sequence of times for which the system switches from the set of active nodes $\G_{k-1}$ to $\G_{k}$.  The notation $T_{\lambda}(u(t))$ refers to the vector $\amatrix{T_{\lambda}(u_1(t)),\ldots,T_{\lambda}(u_N(t))}^T$, and we denote by $F'(t)$ the associated Jacobian matrix with respect to the internal state variables $u$. Because the activation function has the form in \eqref{eq:thresh}, the matrix $F'(t)$ is diagonal with diagonal elements equal to zero for indices in the inactive set $\Gc$ and equal to $f'(u_n)$ for indices $n$ in the active set $\G$. 

Using the chain rule and \eqref{eq:thresh}, we can calculate the derivative with respect to time of the outputs on the active set $\G$:
\begin{equation}
\dot{a}_{\G}(t)= F'_{\G}(t)\dot{u}_{\G}(t)=\textrm{diag}\Big(f'\prn{u_n(t)}\Big)_{\{n\in\G\}}\dot{u}_{\G}(t).
\label{eq:deriv}
\end{equation}
As a consequence, the ODE \eqref{eq:dyn} can be rewritten for nodes in the active set as follows:
\begin{equation} 
\dot{a}_{\G}(t) =F'_{\G}(t)\Prn{-u_{\G}(t)+a_{\G}(t) -\PP a_{\G}(t) + \Phi_{\G}^Ty}. 
\label{eq:active}
\end{equation}
The active nodes follow this ODE until the next switch occurs, changing the sets of active and inactive nodes.  Similarly, we can rewrite the set of differential equations acting on the inactive nodes. Since the output of the activation function for nodes in $\Gc$ is zero, the ODE \eqref{eq:dyn} on the inactive set becomes:
\begin{equation} 
\dot{u}_{\Gc}(t)=-u_{\Gc}(t)-\PcP a_{\G}(t) + \Phi_{\Gc}^Ty .
\label{eq:inactive}
\end{equation}

Separating the system dynamics into two separate characterizations for the active nodes $\G$ and inactive nodes $\Gc$ yields two sets of differential equations that are partially decoupled, and this will be crucial for our subsequent analysis.  This partial decoupling is possible because only the active nodes in the system produce inhibitory signals, coupling their dynamics to the dynamics of each node in the inactive set via the interconnection matrix.  However, because inactive nodes do not inhibit active nodes, the dynamics on the active set are independent of the inactive set.  

\section{Convergence Results}
\label{sec:conv}

In~\cite{rozell_sparse_2008}, the authors show that the LCA network trajectory is non-increasing on the energy surface corresponding to the desired objective function.  While this is necessary behavior for a network that solves an optimization problem, it is not sufficient to actually show that the network converges to a fixed point (or set of points), let alone reaches a minimum of that desired objective function.  Both of these are important guarantees to make before relying on this system as a viable model of neural processing or as an implementation in  engineering applications.


In Section \ref{ssec:conv}, we analyze the general convergence properties of the LCA. In particular, in this section the interconnection matrix $\Phi^T\Phi-I$ may have many negative or zero eigenvalues, which complicates the analysis. The main results of this first section show that the fixed points of the LCA correspond to critical points of the objective \eqref{eq:lp} and the outputs of the network converge from any initial state to the set of fixed points. In the case of isolated critical points of \eqref{eq:lp}, we show that the LCA is globally convergent in the sense that starting from any initial state, the system converges to a fixed point.  Elaborating on these results, Section \ref{ssec:finite} shows that the LCA (even under very general assumptions) is very well-behaved, converging in a finite number of switches (i.e., nodes crossing above or below threshold) and therefore recovering the solution support  (i.e., the set of non-zero elements) in finite time.  Due to the difficulties described above, our approach relies on classic results from nonsmooth analysis, such as subgradients and a generalized chain rule (see Appendix \ref{app:gas}).

\subsection{Global Asymptotic Convergence}
\label{ssec:conv}

To begin, we recall some useful definitions of stability and concepts from Lyapunov theory (see \cite{bacciotti_liapunov_2005} for more background). There exist several notions of stability associated with dynamical systems that describe the evolution of the nodes or their outputs both locally and globally. For a neural network described by a differential equation of the form \begin{equation}
\dot{x}=G(x),\ x\in\field{R}^N,
\label{eq:syst}
\end{equation}
with outputs $z=T_{\lambda}(x)$ defined as in \eqref{eq:thresh}, we say that a constant vector $x^*\in\field{R}^N$ is a \textbf{fixed point} of \eqref{eq:syst} if and only if
\begin{equation}
G(x^*)=0.
\label{eq:fixed}
\end{equation}
On the other hand, the outputs of a system reach a \textbf{critical point} $a^c$ of the objective function $V(\cdot)$ when they satisfy the inclusion:
\begin{equation}
0 \in \partial V(a^c).
\label{eq:crit}
\end{equation}
Note that if $V(\cdot)$ is convex, then the critical points correspond exactly to the global minima of $V$.
We say that the system (\ref{eq:syst}) is (Lyapunov) \textbf{stable} at $x^*$ if for each $\epsilon>0$, there exists $R>0$ such that for all $x_0$ with $\norm{x_0-x^*}<R$, and all the solutions $x(\cdot)$ with initial state $x_0$,
\begin{equation}
\norm{x(t)-x^*}<\epsilon, \ \forall t>0.
\label{eq:stable}
\end{equation}
This is clearly a local notion of stability, guaranteeing that once the trajectory is close to a fixed point, it remains nearby. But this type of stability is insufficient to guarantee global convergence, which means that trajectories approach a fixed point as time goes to infinity. The following notion of stability is slightly stronger, guaranteeing that trajectories approach at least a set of fixed points regardless of the initial state.
\begin{definition}
We say that the outputs of \eqref{eq:syst} are \textbf{globally quasi-convergent} if there exists a set $\mathcal{E}=\{z^*\in\field{R}^N:\ \dot{z}^*=0\}$ such that for all $x_0\in\field{R}^N$, the outputs $z(\cdot)=T_{\lambda}\prn{x(\cdot)}$ with initial state $x_0$ satisfy $\dlims{t\rightarrow+\infty}{z(t)}\in\mathcal{E}$.
\end{definition}

Finally, the strongest and most desirable form of convergence, which guarantees that the nodes converge to a single fixed point, is stated in the next definition.
\begin{definition}
We say that (\ref{eq:syst}) is \textbf{globally convergent}, or equivalently \textbf{globally asymptotically stable}, if there exists a fixed point $x^*$ at which the system is stable, and if for all initial states $x_0\in\field{R}^N$ the solutions $x(\cdot)$ satisfy: $\dlims{t\rightarrow+\infty}{x(t)}=x^*$.
\end{definition}

Typically, global convergence is established through the use of a so-called Lyapunov function $V$. The notation $\dot{V}$ refers to the derivative with respect to time, i.e. $dV/dt$.
\begin{definition}
A function $V:\field{R}^N\mapsto \field{R}$ is a \textbf{weak Lyapunov function} on $\field{R}^N$ if:\\
$(i)$ $V(x)>0$, $\forall x\neq0$;\\
$(ii)$ $V$ is continuous on $\field{R}^N$;\\
$(iii)$ $\dot{V}(x)\leq 0$, $\forall x\in\field{R}^N$; and\\
$(iv)$ $V$ is radially unbounded: $$\dlims{\norm{x}\rightarrow+\infty}{V(x)}=+\infty.$$
\label{def:lyap}
\end{definition}
\noindent Similarly, a function is called a \textbf{\emph{strict Lyapunov function}} if it meets the above conditions, with the exception of having a strict inequality in condition $(iii)$ (i.e., $\dot{V}(x)<0, \forall x\in\field{R}^N$).

\begin{remark} When it is possible to find a weak Lyapunov function for a given dynamical system, the first theorem of Lyapunov \cite{bacciotti_liapunov_2005} guarantees that any fixed point of the system is stable in the sense of \eqref{eq:stable}.   However, to show global convergence of a system, the second theorem of Lyapunov requires a strict Lyapunov function.
\end{remark}
One can check that the objective function in \eqref{eq:lp} is a weak Lyapunov function for the system \eqref{eq:dyn}, thus guaranteeing that the LCA is Lyapunov stable (which improves on the stability result obtained in \cite{rozell_sparse_2008}). However, it is not a strict Lyapunov function.   Indeed, \eqref{eq:lp} only depends on the active nodes, meaning that it could stop decreasing while subthreshold nodes are still evolving.  Continued evolution of the inactive nodes could cause a node in $\Gc$ to become active, thereby causing the objective function to start decreasing again.  As a consequence, condition $(iii)$ in Definition \ref{def:lyap} is not satisfied with strict inequality, and the standard approach of using the objective function as a Lyapunov function is not sufficient to show global convergence of the LCA.  To show global convergence of the system to a fixed point, it is necessary to account for the dynamics on both active and inactive nodes.  

Our main convergence theorem guarantees the global quasi-convergence of the LCA towards a set of critical points of the objective function.  In the case of isolated critical points, this implies global convergence.
\begin{theorem}
The LCA system defined in (\ref{eq:dyn}), with an activation function of the form (\ref{eq:thresh}) satisfying conditions \eqref{eq:cond}:
\begin{enumerate}
\setlength{\itemsep}{0cm}
\setlength{\parskip}{0cm}
 \item has fixed points that are critical points of the objective function defined in (\ref{eq:lp}); 
 \item has globally quasi-convergent outputs; and 
 \item is globally convergent, provided that the critical points of \eqref{eq:lp} are isolated.
\end{enumerate}
\label{th:gas}
\end{theorem}

Note that part 2 of the above theorem only relates to the output variables, while part 3 is a much stronger condition on the entire dynamical system (including subthreshold states).  In the highly relevant case of convex objective functions with unique minima (e.g., BPDN), part 3 of the above theorem applies directly and we have the strongest possible notion of convergence.  In fact, for most dictionaries $\Phi$ (e.g., random Gaussian), the minimum of the $\ell_1$-minimization objective function \eqref{eq:l1} is unique \cite{fuchs_sparse_2004}. Recent results on subgradient dynamical systems \cite{bolte_lojasiewicz_2007,forti_convergence_2006} lead us to believe that part 3 could still hold in the case where the fixed points of the system are not isolated and there exists a subspace of solutions to \eqref{eq:lp}, but this conjecture is beyond the scope of this work.  Section \ref{ssec:exConv} illustrates the convergence behavior of the LCA in simulation.  The proof of this theorem is in Appendix~\ref{app:gas}, and relies on generalized notions of a subgradient due to the nonsmooth nature of the objective.




\subsection{Convergence in a Finite Number of Switches}
\label{ssec:finite}

In the theorem below, we strengthen the previous result to establish that under some mild conditions, the support of the final solution is reached with a finite number of switches.  To prove this result, it is sufficient to assume that no node in the solution lies exactly on the threshold $\lambda$. This assumption precludes unwanted infinite oscillation behavior on the boundaries, known as Zeno behavior. In other words, we assume that there exists a margin ($r>0$) above and below the threshold which contains no node in $u^*$.  One would expect this condition to hold with near certainty for any signal that was not pathologically constructed.

\begin{theorem}
\label{th:finite}
If the system \eqref{eq:dyn} converges to a fixed point $u^*$ such that there exists $r>0$:
\begin{equation*}
\sarray{\forall n\in\G_*,\ \ \abs{u^*_n}\geq \lambda + r\\
\forall n\in\G_*^c,\ \ \abs{u^*_n}\leq \lambda - r},
\end{equation*}
then the system converges after a finite number of switches.
\end{theorem}
This implies that the neural network recovers the support of the solution $a^*$ in finite time.  Section \ref{ssec:exFinite} explores the number of switches during the convergence of the LCA in simulation.

\begin{IEEEproof}
Let $\G_*$ be the set of active nodes in $u^*$.
By contradiction, assume that the sequence of switching times $\left\{t_k\right\}_{k\in\field{N}}$ is infinite. Since the LCA converges to $u^*$, we have:
$$u(t_k)\underset{k\rightarrow+\infty}{\longrightarrow}u^* .$$
As a consequence, for $r>0$, there exists $K\in\field{N}$ such that $\forall k\geq K$, $\eltwo{u(t_k)-u^*}<r$.
To begin, we show that for all $k\geq K$, the state variables $u(t_k)$ are in the subsystem $\G_*$. For all $k\geq K$, we have two cases:
\begin{itemize}
\item Nodes that are above threshold in $u^*$ are above threshold in $u(t_k)$. Indeed, $\forall n\in\G_*$, we have:
$$r>\abs{u_n(t_k)-u^*_n} \geq\abs{u^*_n}-\abs{u_n(t_k)}\geq\lambda+r-\abs{u_n(t_k)}$$
$$\Rightarrow \ \abs{u_n(t_k)}>\lambda . $$
Moreover, nodes are active with the correct sign, otherwise, we would have:
$$r>\abs{u_n(t_k)-u^*_n} =\abs{u_n(t_k)}+\abs{u^*_n} >\lambda+\lambda+r$$
$$\Rightarrow \ 0>\lambda,$$
which is a contradiction.
\item Nodes that are below threshold in $u^*$ are below threshold in $u(t_k)$. Indeed, $\forall n\in{\G_*^c}$, we have:
$$\abs{u_n(t_k)}-\lambda \leq\abs{u_n(t_k)}-\abs{u^*_n}-r \leq\abs{u_n(t_k)-u^*_n}-r< 0$$
$$\Rightarrow \ \abs{u_n(t_k)}<\lambda .$$
\end{itemize}
As a consequence, for all $k\geq K$, \ \ $\G_k=\G_*$.   However, $\G_k$ and $\G_{k+1}$ must be different to define the switching time $t_{k+1}$ and we reach a contradiction. This proves that after a finite number of switches $K$, there cannot be any switching out of subsystem $\G_*$.
\end{IEEEproof}

\section{Exponential Convergence Rate}
\label{sec:rate}
With Theorem~\ref{th:gas} showing that the LCA is globally convergent, the most pressing issue remaining is to determine the convergence rate of the network.  Such a bound will be especially important for implementations, which must guarantee solution times.  In this section, we show that under some additional conditions on the problem's specifics (given in~\eqref{eq:cond5}), the LCA network converges exponentially fast to a unique fixed point $u^*$. We also give an analytic characterization of the convergence speed.\footnote{We reintroduce the time constant $\tau$ in this discussion of the LCA, since it appears in the expression for the convergence speed.}

To begin this section, we recall the definition of exponential convergence:

\begin{definition}
The dynamical system in \eqref{eq:syst} is \textbf{exponentially convergent} to the solution $x^*$ if there exists a constant $c>0$ such that for any initial point $x(0)$, there exists a constant $\kappa_0>0$ (which may depend on $x(0)$) for which the trajectory $x(t)$ of the system satisfies:
$$\norm{x(t)-x^*} \leq \kappa_0 e^{-ct}, \qquad \forall t\geq0.$$
The constant $c$ is referred to as \textbf{convergence speed} of the system.
\end{definition}

In order to state the main theorem of this section, we define the two following quantities.
%
The first constant is denoted by $\alpha$ and provides a bound on the derivative of $f(\cdot)$ in \eqref{eq:thresh}:
\begin{equation}
\forall t\geq0, \forall n=1,\ldots,N \qquad \abs{f'\prn{u_n(t)}}\leq\alpha.
\label{eq:cond4}
\end{equation}
Note that the constant $\alpha$ is always well defined since the trajectories $u_n(t)$ are bounded and the function $f(\cdot)$ is continuous. The second constant, denoted by $\delta$, is the smallest positive constant such that for any active set $\G$ visited by the algorithm and any vector $x$ in $\field{R}^N$ with active set $\widetilde{\G}=\G\cup\G_*$, where $\G_*$ is the active set of the solution to \eqref{eq:lp}, we have:
\begin{equation}
\prn{1-\delta} \eltwo{x}^2 \leq \eltwo{\Phi x}^2 \leq \prn{1+\delta} \eltwo{x}^2.
\label{eq:delta}
\end{equation}
The constant $\delta$  depends on the singular values of the matrix $\Phi_{\widetilde{\G}}$ and on the sequence of active sets visited by the system. It may not be well defined for any matrix $\Phi$ or any input $y$. However, in many interesting cases in CS, the constant $\delta$ is close to $0$ and the dictionary elements are almost orthogonal for any small enough active set~\cite{candes_stable_2006}. The following theorem shows that this constant directly relates to the convergence speed of the neural network.

\begin{theorem}
\label{th:rate}
Under conditions \eqref{eq:cond} on the activation function in \eqref{eq:thresh}, and provided that the constants $\alpha$ and $\delta$, defined in \eqref{eq:cond4} and \eqref{eq:delta} satisfy 
\begin{equation}
\alpha \delta < 1,
\label{eq:cond5}
\end{equation} 
the LCA network defined in \eqref{eq:dyn} is globally exponentially convergent to a unique equilibrium with convergence speed
$$c = \dfrac{1-\alpha\delta}{\tau}.$$
\end{theorem}

If condition \eqref{eq:cond5} is satisfied, the expression given for the convergence speed is positive and thus meaningful. 
It depends on the eigenvalues of the matrix $\Phi_{\G}^T\Phi_{\G}$, which vary with the active set $\G$. A careful analysis of the sequence of active sets visited by the network is required to obtain a good estimate of $\delta$. Since such a study is application dependent, we do not address this question here.
Note that in the very interesting case of the soft-threshold function, $\alpha=1$ and so, condition \eqref{eq:cond5} reduces to $\delta<1$. The time constant $\tau$ of the physical solver implementing the LCA neural network appears in the expression of the speed of convergence. Lowering this time constant means the system will converge faster.  Analog systems can have smaller time constants than their digital counterparts that scale better with the problem size. Section \ref{ssec:exRate} explores the convergence rate bounds for the LCA in simulation.

To establish the expression of the convergence speed, we use the following energy function:
\begin{equation}
E(t)=\dfrac{1}{2}\eltwo{\widetilde{u}(t)}^2,
\label{eq:energy2}
\end{equation}
where we again redefine the output and state variables in terms of the distance from any arbitrary fixed point $u^*$ of \eqref{eq:dyn}
\begin{equation}
\begin{split}
\widetilde{u}_n(t) & = u_n(t)-u_n^*,\\
\widetilde{a}_n(t) & = a_n(t) - a_n^* = T_{\lambda}(\widetilde{u}_n(t)+u_n^*)-T_{\lambda}(u_n^*).
\end{split}
\label{eq:var2}
\end{equation}
We will show that the function \eqref{eq:energy2} converges exponentially fast to zero by first analyzing it on the set of indices $\widetilde{\G}$. From this, we obtain a convergence result for the outputs. We then use this result to prove convergence of the entire state vector to a unique fixed point $u^*$. The properties presented in the following Lemma (see Appendix \ref{app:pties} for proof) are useful to prove the main result of this section.

\begin{lemma}
\label{lem:pties}
Assume that the activation function \eqref{eq:thresh} satisfies the conditions \eqref{eq:cond}. Then, the set of variables $\widetilde{u}$ and $\widetilde{a}$ defined in \eqref{eq:var2} satisfy the properties:
\begin{enumerate}[(i)]
\item $\sign{\widetilde{a}_n} = \sign{\widetilde{u}_n}.$
\item $\abs{\widetilde{a}_n} \leq \alpha\abs{\widetilde{u}_n}.$
\item $\widetilde{a}_{\mathcal{T}}^T\widetilde{a}_{\mathcal{T}} \leq \alpha\widetilde{u}_{\mathcal{T}}^T\widetilde{a}_{\mathcal{T}} \leq \alpha^2\widetilde{u}_{\mathcal{T}}^T\widetilde{u}_{\mathcal{T}}$ \ \ for any $\mathcal{T}$\\
(in particular for $\mathcal{T}=\widetilde{\G}$).
\item $\dsum{n=1}{N}{\dint{0}{\widetilde{u}_n(t)}{g_n(s) ds}}\leq\widetilde{u}^T\widetilde{a}$.
\end{enumerate}
\end{lemma}

Armed with this lemma, we can now give the proof of the main theorem.

\begin{IEEEproof}[Proof of Theorem \ref{th:rate}]
Using the fact that $u^*$ is a fixed point, we set $\dot{u}^*(t)=0$ in \eqref{eq:dyn} and rewrite the dynamics in terms of the new variables:
\begin{equation}
\tau\dot{\widetilde{u}}(t) = -\widetilde{u}(t) - \prn{\Phi^T\Phi - I}\widetilde{a}(t).
\label{eq:tildesol}
\end{equation}
We start by showing that the partial energy function $E_{\widetilde{\G}}(t) = \dfrac{1}{2}\eltwo{\widetilde{u}_{\widetilde{\G}}(t)}^2$ converges exponentially fast, where $\widetilde{\G}$ is the support of $\widetilde{a}(t)$. Using the chain rule, the time derivative of $E_{\widetilde{\G}}(t)$ along the network trajectory is
\begin{align*}
\tau \dot{E}_{\widetilde{\G}}(t) & = \tau \dot{\widetilde{u}}_{\widetilde{\G}}^T(t) \widetilde{u}_{\widetilde{\G}}(t) \\
& = - \widetilde{u}_{\widetilde{\G}}^T(t)\prn{\widetilde{u}_{\widetilde{\G}}(t) + \prn{\Phi^T_{\widetilde{\G}}\Phi_{\widetilde{\G}} -I_{\widetilde{\G}}} \widetilde{a}_{\widetilde{\G}}(t)} \\
& = -\eltwo{\widetilde{u}_{\widetilde{\G}}(t)}^2  - \widetilde{u}_{\widetilde{\G}}^T(t) \prn{\Phi^T_{\widetilde{\G}}\Phi_{\widetilde{\G}}-I_{\widetilde{\G}}} \widetilde{a}_{\widetilde{\G}}(t).
\end{align*}
Since $\widetilde{a}_{\widetilde{\Gamma}}$ is supported on ${\widetilde{\Gamma}}$, assumption \eqref{eq:delta} implies that the eigenvalues of $  \Phi_{\widetilde{\Gamma}}^T\Phi_{\widetilde{\Gamma}}$ lie between $(1-\delta)$ and $(1+\delta)$ and so:
\begin{align*}
  \left\| \left( \Phi_{\widetilde{\Gamma}}^T\Phi_{\widetilde{\Gamma}} - I_{\widetilde{\Gamma}} \right) \widetilde{a}_{\widetilde{\Gamma}}  \right\|_2  \leq  \left\| \Phi_{\widetilde{\Gamma}}^T\Phi_{\widetilde{\Gamma}} - I_{\widetilde{\Gamma}} \right\|_2 \| \widetilde{a}_{\widetilde{\Gamma}} \|_2 & \\
   \leq \max \left\lbrace (1+\delta) -1 , \; 1- (1-\delta) \right\rbrace \| \widetilde{a}_{\widetilde{\Gamma}} \|_2 &\\
  = \delta \| \widetilde{a}_{\widetilde{\Gamma}} \|_2. &
\end{align*}
Finally, property $(iii)$ of Lemma \ref{lem:pties} states that for any set $\mathcal{T}$, $\eltwo{\widetilde{a}_{\mathcal{T}}}^2 \leq \alpha^2 \eltwo{\widetilde{u}_{\mathcal{T}}}^2$. Using the Cauchy-Schwartz inequality and putting everything together, we obtain:
\begin{align*}
\abs{\widetilde{u}_{\widetilde{\G}}^T \prn{\Phi^T_{\widetilde{\G}}\Phi_{\widetilde{\G}}-I_{\widetilde{\G}}} \widetilde{a}_{\widetilde{\G}} } & \leq \eltwo{\widetilde{u}_{\widetilde{\G}}} \eltwo{ \prn{\Phi^T_{\widetilde{\G}}\Phi_{\widetilde{\G}}-I_{\widetilde{\G}}} \widetilde{a}_{\widetilde{\G}} } \\
 & \leq \eltwo{\widetilde{u}_{\widetilde{\G}}} \delta \eltwo{\widetilde{a}_{\widetilde{\G}}} \\
 & \leq \alpha \delta \eltwo{\widetilde{u}_{\widetilde{\G}}}^2.
\end{align*}
As a consequence, we obtain that the time derivative of the partial energy function satisfies
\begin{align*}
\tau \dot{E}_{\widetilde{\G}}(t) & \leq  -\eltwo{\widetilde{u}_{\widetilde{\G}}(t)}^2  +  \alpha \delta \eltwo{\widetilde{u}_{\widetilde{\G}}(t)}^2 \\
& \leq - 2\prn{1-\alpha\delta} E_{\widetilde{\G}}(t).
\end{align*}
Using Gronwall's inequality on the interval $\Prn{t_k,t_{k+1}}$ where $\widetilde{\G}$ is constant yields
$$E_{\widetilde{\G}}(t) = \dfrac{1}{2}\eltwo{\widetilde{u}_{\widetilde{\G}}(t)}^2 \leq  \dfrac{1}{2}\eltwo{\widetilde{u}_{\widetilde{\G}}(t_k)}^2e^{-2(1-\alpha\delta)(t-t_k)/\tau}.$$
Since $\eltwo{\widetilde{a}(t)} \leq \alpha \eltwo{\widetilde{u}_{\widetilde{\G}}(t)}$, $\forall t\in\Prn{t_k,t_{k+1}}$:
\begin{equation}
\eltwo{\widetilde{a}(t)} \leq \alpha \eltwo{\widetilde{u}_{\widetilde{\G}}(t_k)}e^{-(1-\alpha\delta)(t-t_k)/\tau}.  \label{eq:outputdecay}
\end{equation}
%
%

Using this result on the output, we now prove that the state $u(t)$ converges exponentially fast. For this we write the solution to \eqref{eq:tildesol} as follows:
\begin{align*}
\widetilde{u}(t) & = e^{-(t-t_k)/\tau}\widetilde{u}(t_k) \\
& \qquad + e^{-(t-t_k)/\tau}\dint{t_k}{t}{ e^{(\nu-t_k)/\tau} \prn{I-\Phi^T\Phi}\widetilde{a}(\nu) d\nu}. 
\end{align*}
Denoting by $h(t)$ the second term in the right-hand side, and plugging in \eqref{eq:outputdecay}, we bound the norm of $h(t)$ by
\begin{align*}
 \eltwo{h(t)} & \leq  e^{-(t-t_k)/\tau}\dint{t_k}{t}{e^{(\nu-t_k)/\tau}\eltwo{\prn{\Phi^T\Phi-I} \widetilde{a}(\nu)} d\nu} \\
&  \leq e^{-(t-t_k)/\tau}\dint{t_k}{t}{ \underbrace{\eltwo{\Phi^T\Phi-I}}_{=C_1} e^{(\nu-t_k)/\tau}\eltwo{\widetilde{a}_{\widetilde{\G}}(\nu)}   d\nu} \\
&  \leq e^{-(t-t_k)/\tau}\dint{t_k}{t}{ C_1\alpha\eltwo{\widetilde{u}_{\widetilde{\G}}(t_k)} e^{\alpha\delta(\nu-t_k)/\tau}  d\nu} \\
&  =  \dfrac{C_1\tau}{\delta} \eltwo{\widetilde{u}_{\widetilde{\G}}(t_k)} e^{-(t-t_k)/\tau} \Prn{e^{\alpha\delta (t-t_k)/\tau} -1} \\
%
%
&  \leq C_2\eltwo{\widetilde{u}_{\widetilde{\G}}(t_k)} ~e^{-(1-\alpha\delta)(t-t_k)/\tau}\\
&  \leq C_2\eltwo{\widetilde{u}(t_k)} ~e^{-(1-\alpha\delta)(t-t_k)/\tau},
\end{align*}
where $C_2 = \prn{\eltwo{\Phi^T\Phi-I}\tau / \delta}$. We plug this bound in the expression for $\widetilde{u}(t)$ to get:
\begin{align*}
\eltwo{\widetilde{u}(t)} & = \eltwo{ e^{-(t-t_k)/\tau}\widetilde{u}(t_k) + h(t) }\\
& \leq \eltwo{\widetilde{u}(t_k)}e^{-(t-t_k)/\tau} + \eltwo{h(t)} \\
& \leq \eltwo{\widetilde{u}(t_k)} e^{-(t-t_k)/\tau} \\
& \qquad\qquad\qquad + C_2 \eltwo{\widetilde{u}(t_k)} e^{-(1-\alpha\delta)(t-t_k)/\tau} \\
 & \leq \prn{1+ C_2} \eltwo{\widetilde{u}(t_k)} e^{-(1-\alpha\delta)(t-t_k)/\tau} \\
 & = C_3 \eltwo{\widetilde{u}(t_k)} e^{-(1-\alpha\delta)(t-t_k)/\tau},
\end{align*}
where $C_3 = 1+C_2$. Since $ \eltwo{\widetilde{u}(t)} $ is continuous for all time $t$, it is easy to show (by induction on $t_k$) that
\begin{equation}
 \eltwo{\widetilde{u}(t)} \leq e^{-(1-\alpha\delta)t/\tau} C_3 \eltwo{ \widetilde{u}(0)}. \label{eq:decay}
\end{equation}
This last inequality shows that the state variable converges exponentially fast to a unique fixed point $u^*$ with convergence speed $(1-\alpha\delta)/\tau$.
\end{IEEEproof}


\section{Simulations}
\label{sec:sim}

In this section, we illustrate the previous theoretical results in simulation. Each plot is based on the following canonical sparse approximation problem.
We generate a ``true'' sparse vector $a_0\in\field{R}^N$, with $N=512$ and $s=5$ non-zero entries.  We select the locations of the nonzeros uniformly at random and draw their amplitudes from a normal gaussian distribution (then normalizing to have unit norm). We choose a dictionary $\Phi$ as a union of the canonical basis and a sinusoidal basis having dimensions $M\times N$ with $M=256$. The vector of measurements is $y=\Phi a_0 + \eta \in\field{R}^M$, where $\eta$ is Gaussian random noise with standard deviation $0.0062$. We use an LCA with a soft-threshold activation function, with a threshold set to $\lambda=0.025$ and $u(0)=0$.  We simulate the LCA dynamical equations \eqref{eq:dyn} through a discrete approximation in Matlab with a step size of $0.001$ and a solver time constant chosen to be equal to $\tau=0.01$.

\subsection{Convergence Results}
\label{ssec:exConv}
\begin{figure}
\centering
\includegraphics[width=3in]{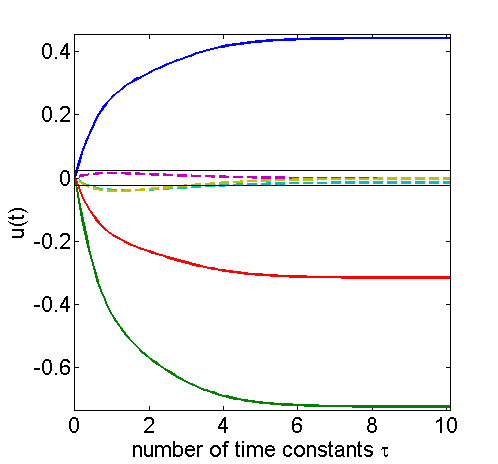}
\caption{Plot of the evolution with respect to time of several LCA nodes $u_k(t)$. The plain lines correspond to nodes that are active in the final solution and the dashed lines correspond to nodes that are inactive in the final solution.}
\label{fig:conv2}
\end{figure}

\begin{figure}
\centering
\includegraphics[width=3in]{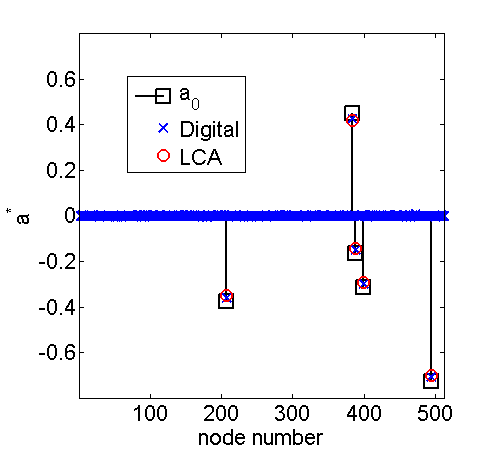}
\caption{Output $a^*$ of the LCA after convergence. Only non-zero elements are plotted. The fixed point reached by the system is very close to the initial sparse vector used to create the measurement vector (it cannot be exact due to noise).  The solution is also very close to a standard digital solver~\cite{kim_interior-point_2007} run using the same inputs.  Note that the LCA produces many coefficients that are exactly zero (therefore not plotted).}
\label{fig:conv1}
\end{figure}

From Theorem \ref{th:gas}, the LCA should converge and recover the solution to the sparse approximation problem \eqref{eq:lp}, which has a unique minimizer. Since the signal $a_0$ to recover is sparse, the outputs of the neural network are expected to converge to a solution close to the initial signal $a_0$. \figurename~\ref{fig:conv2} shows the evolution of a few nodes $u_n(t)$ selected at random. We see that both active and inactive nodes converge relatively quickly.  \figurename~\ref{fig:conv1} shows the fixed point reached by the LCA system and compares it to the initial signal $a_0$ and to the solution obtained using a standard digital solver for the same optimization program \cite{kim_interior-point_2007}. The solution reached by the network possesses exactly 5 non-zero entries that correspond to the non-zero entries in $a_0$.  The recovered amplitudes are very close to the initial amplitudes (it cannot be exact due to the added measurement noise) and to the ones produced by the reference digital solver. However, the LCA produces a sparse vector while the digital solver returns many small but non-zero entries that would have to be removed by postprocessing. 

\begin{figure}
\centering
\includegraphics[width=3in]{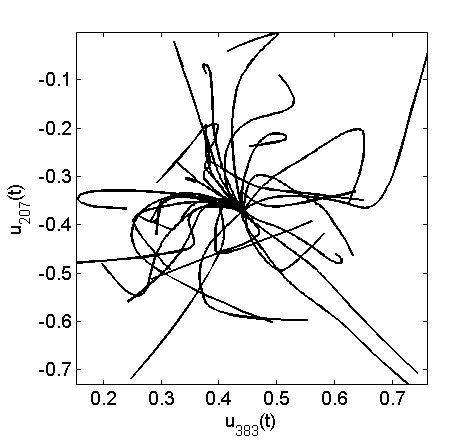}
\caption{Trajectories of $u(t)$ in the plane defined by two nodes chosen randomly from the active set.  Trajectories are shown for 30 random initial states.}
\label{fig:traj}
\end{figure}

To illustrate the global convergence behavior, we also ran the LCA for 30 randomly generated initial points. We selected two nodes from the final active set and plotted the trajectories in the space defined by those two nodes. \figurename~\ref{fig:traj} clearly shows that the solution is attractive for any of those initial points.

\begin{figure}
\centering
\includegraphics[width=3.4in]{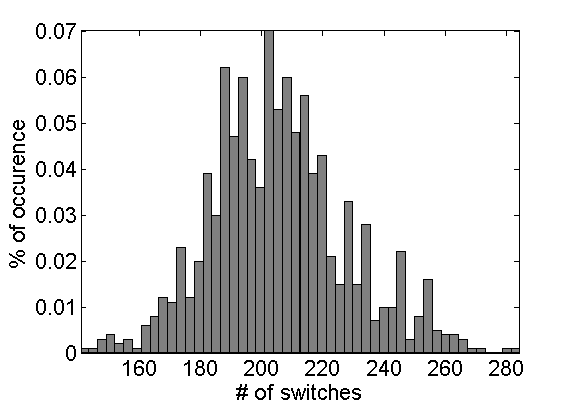}
\caption{Histogram (in percentage) of the number of switches the LCA requires before convergence over 1000 trials.}
\label{fig:hist}
\end{figure}

\label{ssec:exFinite}
To illustrate the number of switches used by the system (see Theorem \ref{th:finite}), we generate 1000 sparse vectors $a_0$ and measurements $y$ and simulate the LCA dynamics. \figurename~\ref{fig:hist} shows a histogram of the number of switches needed for the system to converge. The figure illustrates that the number of switches before convergence of the neural network is finite and of the order of the dictionary size. 
This illustrates that this solver takes an efficient path towards the solution.

\subsection{Convergence Rate Results}
\label{ssec:exRate}

\begin{figure}
\centering
\includegraphics[width=3in]{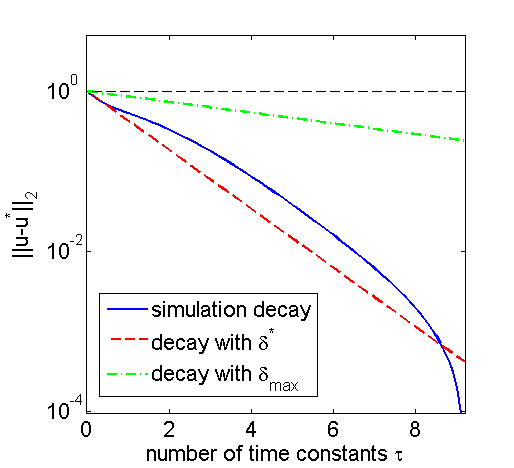}
\caption{Convergence behavior of $\norm{u(t)-u^*}_2$. The dashed line shows the theoretical decay in \eqref{eq:decay} with $\delta^*$ computed by using the final solution support. The dash-dot line shows the theoretical decay with $\delta_{\max}$ computed on the largest support visited.  While both estimates are showing theoretically correct behavior, the estimate of the rate based on $\delta^*$ is more empirically accurate than the conservative estimate based on $\delta_{\max}$.}
\label{fig:rate}
\end{figure}

To illustrate the convergence rate result in Theorem \ref{th:rate}, it is necessary to find an expression for the convergence speed $(1-\delta)/\tau$ that appears in the exponential term in \eqref{eq:decay} (here $\alpha=1$). This term bounds the error term: $$\norm{u(t)-u^*}_2,$$ which is normalized to have initial value of $1$ at $t=0$ and plotted using a log-scale on \figurename~\ref{fig:rate}. Note that the constant $\delta$ defined in \eqref{eq:delta} depends on the sequence of active sets $\G_k$ visited by the system. However, it is very difficult to predict for a given input signal what sequence of active sets the algorithm is going to visit. To estimate this upper bound, we compute the constant $\delta$ using the matrix $\Phi_{\G_*}$ composed of the dictionary elements that are active in the final solution. The corresponding upper bound on the decay $e^{-(1-\delta^*)t/\tau }$ is plotted in \figurename~\ref{fig:rate}. During the transient phase, the number of active nodes is actually larger and thus, we expect $1-\delta$ to be smaller than this estimate and the nodes to converge slower at some point during the transient phase. As a consequence, we keep track of the largest support visited by the network and compute the corresponding $\delta$. This second upper bound $e^{-(1-\delta_{\max})t/\tau}$ on the convergence rate is plotted on the same figure. As expected, the theoretical decay computed with the maximum support visited is an upper bound for the convergence speed. However, it can be seen that this estimate is very pessimistic and that the bound computed with the final active set is a better estimate for the experimental decay. This simulation illustrates that the  theoretical exponential convergence appears to capture the essential system behavior.

\section{Conclusion}
\label{sec:ccl}

This paper presents a mathematical analysis of the convergence properties and convergence rate of the LCA, a neural network designed specifically for the challenging  sparse approximation problem. Despite a nonsmooth activation function and possibly singular interconnection matrix that prevent the application of existing analysis approaches, we have shown that the system is globally convergent to the optimal solution.  In addition, under some mild assumptions on the solution, we have shown that the trajectories follow a reasonable path and reach the final active set in finite time. Finally, under slightly stronger assumptions on the problem specifics (applicable at least in CS recovery problems), we established that the LCA is exponentially convergent and with a convergence rate that depends on problem-specific parameters. 

This collection of results and analysis leads us to conclude that performance guarantees can be made for the LCA system that make it plausible for implementation in engineering applications and as a model of biological information processing.  Indeed, providing such guarantees makes it easier to justify the expense associated with developing analog VLSI implementations, which could eventually result in significant improvements to the speed and power consumption necessary for real-time signal processing applications.  Our future work will concentrate on finding reasonable estimates of the theoretical convergence speed (especially in well-studied special cases such as CS recovery), and further characterizations of the LCA system dynamics that may open up new applications of this system for time-varying input signals.

\appendices
\section{Proof of Theorem \ref{th:gas}}
\label{app:gas}

Building on the earlier definition of a subgradient to give a notion that non-differentiable functions can still be well-behaved,  a function $g:X\rightarrow \field{R}$ (where $X$ is a Banach space) is said to be \emph{regular}~\cite[Def. 2.3.4]{clarke_optimization_1987} at $x$ in $X$ if
\begin{enumerate}[(i)]
\item For all $v\in X$, the usual one-sided  directional derivative
$$g'(x;v)=\lim_{t\downarrow 0} \dfrac{g(x+tv)-g(x)}{t} \qquad \mathrm{exists.}$$
\item For all $v \in X$, $g'(x;v)=g^{\circ}(x;v)$, where $g^{\circ}(x;v)$ is the generalized directional derivative
$$g^{\circ}(x;v)=\limsup_{\substack{y\rightarrow x \\ t\downarrow 0}} \dfrac{g(y+tv)-g(y)}{t}.$$
\end{enumerate}
It can easily been seen that, since the function $C(\cdot)$ is defined from $\field{R}$ to $\field{R}$, is differentiable on $\field{R}\backslash\{0\}$ and the function $T_{\lambda}(\cdot)$ is continuous on all $\field{R}$, then $C(\cdot)$ admits left and right derivatives and is clearly regular on $\field{R}$.  This implies that $V(\cdot)$ is regular on $\field{R}^N$, and by \cite[Prop. 2.3.3]{clarke_optimization_1987} we have that:
\begin{equation}
\begin{split}
	\partial{V(a(t))} & = -\Phi^Ty + \Phi^T\Phi a(t) + \lambda \partial{C(a(t))},
\end{split}
\label{eq:grad}
\end{equation}
where $\partial{C(a(t))}=\amatrix{\partial{C(a_1(t))},\ldots,\partial{C(a_N(t))}}^T$.
We also recall the following result \cite[Th. 2.3.9 (iii)]{clarke_optimization_1987}, which is a generalization of the chain rule for regular functions.
\begin{lemma}
Suppose that $V(a): \field{R}^N \rightarrow \field{R}$ is regular in  $\field{R}^N$ and that $a(t): \left[0, +\infty \right) \rightarrow \field{R}^N$ is strictly differentiable on $\left[0, +\infty \right)$. Then, $V(a(t))$ is regular on $\field{R}^N$, and we have
\begin{equation}
\dfrac{d}{dt} V(a(t)) = \zeta^T\dot{a}(t) \qquad \forall \zeta \in \partial V(a(t)). 
\label{eq:chain}
\end{equation}
\end{lemma}
Note that from this theorem, since $V(a(t))$ is regular, we can choose any element in $\partial V(a(t))$ to compute the time derivative of $V(\cdot)$ along the trajectories of the neural network. Armed with these tools, we proceed with the proof of Theorem~\ref{th:gas}.

\begin{IEEEproof}[Proof of Theorem \ref{th:gas}]
Beginning with part 1 of the theorem, we first show that any fixed point of system \eqref{eq:dyn} is a critical point of the objective in (\ref{eq:lp}).
From \eqref{eq:fixed}, any fixed point $u^*$ of \eqref{eq:dyn} satisfies the relationship $\dot{u}(t)=0$. Let $\G_*$ be the active set at the fixed point $u^*$. Eq. \eqref{eq:active} and \eqref{eq:inactive} yield:
\begin{align}
 & -u^*_{\G_*} + f(u^*_{\G_*}) - \PP a^*_{\G_*} + \Phi_{\G_*}^Ty = 0 \label{eq:fixa} \\
 & -u^*_{\G_*^c}-\Phi_{\G_*^c}^T\Phi_{\G_*} a^*_{\G_*} + \Phi_{\G_*^c}^Ty = 0.
\label{eq:fixi}
\end{align}
According to \eqref{eq:crit}, and using \eqref{eq:grad}, a point $a^*$ is a critical point of $V(\cdot)$ if and only if:
\begin{equation}
\Phi^Ty - \Phi^T\Phi a^* \in \lambda \partial{C(a^*)}.
\label{eq:nsc}
\end{equation}
For the nodes in the active set $n\in\G_*$, $C(a_n)$ admits a usual gradient, as defined in \eqref{eq:corres} and thus, \eqref{eq:nsc} yields
$$ \Phi^T_{\G_*}y - \Phi_{\G_*}^T\Phi_{\G_*}a^*_{\G_*} = u^*_{\G_*} - f(u^*_{\G_*}). $$
This is the same as condition \eqref{eq:fixa}. For nodes in the inactive set $\G_*^c$, we need to determine the subgradient of $C(\cdot)$ at zero. To do so, note that using \eqref{eq:corres} and the continuity of the function $f(\cdot)$ at $\lambda$ (condition \eqref{eq:cond1}), we have:
$$ \lim_{a\downarrow 0} \lambda\nabla C(a) = \lim_{u \downarrow \lambda} u - f(u) = \lambda. $$
As a consequence, we find that:
$$ \lim_{a\downarrow 0} \nabla C(a) = 1 \qquad \mathrm{and} \qquad \lim_{a\uparrow 0} \nabla C(a) = -1 .$$
From the definition of the subgradient, these limits show that:
$\partial{C(0)}=\Prn{-1,1}.$
Using this, the condition in~\eqref{eq:nsc} restricted to the set of inactive nodes is equivalent to:
$$ \Phi^T_{\G_*^c}y - \Phi_{\G_*^c}^T\Phi_{\G_*} a^*_{\G_*} \in [-\lambda, \lambda] .$$
We immediately see that this condition is the same as \eqref{eq:fixi}, since by definition of the inactive nodes $u^*_{\G_*^c}\in[-\lambda,\lambda]$. This shows that the fixed points $a^*$ coincide with the critical points of the objective function \eqref{eq:lp}.

Moving on to establishing the convergence result in part 2 of the theorem, we first note that from our analysis above, the same reasoning leads to the conclusion that:
\begin{equation}
-\dot{u}(t) = u(t) - f(u(t))  + \Phi^T \Phi \ a(t) - \Phi^Ty \in \partial{V(a(t))}.
\label{eq:parV}
\end{equation}
Using the chain rule \eqref{eq:chain} with $\zeta=-\dot{u}(t) \in \partial{V(a(t))}$ from \eqref{eq:parV} and using the expression in \eqref{eq:deriv}, we get the relationship:
\begin{align} 
 \dff{V(a(t))}{t} & = -\dot{u}(t)^T\dot{a}(t) \nonumber \\
 & = -\dot{u}(t)^TF'(t)\dot{u}(t) \nonumber \\
 & = -\sum\limits_{n\in\G} f'(u_n(t)) \abs{\dot{u}_{n}(t)}^2.
\label{eq:dva}
\end{align}
This expression is valid for all time $t\geq0$. In addition, because the function $f(\cdot)$ satisfies \eqref{eq:cond2}, $f'(u_n(t))>0$, and thus 
\begin{align*}
 & \dff{V(a(t))}{t}\leq0 & & \text{for all } t\geq0, \mbox{ and}\\
 & \dff{V(a(t))}{t}<0 & & \text{for all } t\geq0 \text{ such that } \eltwo{\dot{u}_{\G}(t)}\neq0.
\end{align*}
This means that $V(a(t))$ is non-increasing for all $t\geq0$. Since $V(a(t))$ is continuous, bounded below by zero, and non-increasing, $V(a(t))$ converges to a constant value $V^*$, and its time derivative  $\dot{V}(a(t))$ tends to zero as $t\to\infty$. Using equations (\ref{eq:dva}) and \eqref{eq:cond2}, we conclude that:
$\dlims{t\rightarrow+\infty}{\eltwo{\dot{u}_{\G}(t)}}=0.$
As a consequence, we also have
$\dlims{t\rightarrow+\infty}{\eltwo{\dot{a}(t)}}=0$,
so the outputs converge to the set
$E = \left\{a: \mathrm{\ s.t.\ } \dot{a}(t)=0\right\},$
and the LCA outputs are quasi-convergent.

Moving on to part 3 of the theorem, we assume that the critical points of \eqref{eq:lp} are isolated. We need to show that both active and inactive nodes converge to a single fixed point. From part 2, we know that since the nodes converge to $E$, after some time, the active nodes will be within a ball of radius $R$ around one element $a^* \in E$. However, since critical points of \eqref{eq:lp} are isolated, there exists a ball $\mathcal{B}_{\epsilon}(a^*)$ of radius $\epsilon>0$ around $a^*$ that does not contain any critical point: 
$$a^*\in E \qquad \text{and} \qquad \forall a\in\mathcal{B}_{\epsilon}(a^*), a\neq a^* \Rightarrow a\notin E.$$
Since the system is stable, we know that once the trajectory gets close enough (within a ball of radius $R$) to one element in $E$, it cannot leave a ball of radius $\epsilon$ around this fixed point (see \eqref{eq:stable}). As a consequence, the outputs remain within the ball $\mathcal{B}_{\epsilon}(a^*)$, which contains only the fixed point $a^*$. This proves that the active nodes converge to the point $a^*$ in $E$:
\begin{equation}
\dlims{t\rightarrow+\infty}{a(t)}=a^*.
\label{eq:activeConv}
\end{equation}
This implies that the ODE~\eqref{eq:dyn} can be written in terms of the distance $\widetilde{a}(t) = a(t)-a^*$ of the outputs from the solution:
$$\dot{u}(t)=-u(t)-\Phi^T\Phi a^* + \Phi^Ty + a^* - \Phi^T\Phi \widetilde{a}(t) +\widetilde{a}(t).$$
We let
$$u^* = -\Phi^T\Phi a^* + \Phi^Ty + a^*,$$
and rewrite the ODE \eqref{eq:dyn} as:
$$\dot{u}(t)=-u(t)+u^* - \prn{\Phi^T\Phi-I} \widetilde{a}(t)$$
Solving this ODE for all $t\geq0$ yields:
$$u(t) = u^*+e^{-t}\prn{u(0)-u^*} + e^{-t}\dint{0}{t}{e^s\prn{\Phi^T\Phi-I} \widetilde{a}(s) ds}.$$
While it is difficult to say anything directly about the trajectory of the system, it is helpful to consider a surrogate trajectory that is a straight line in the state-space: $u^*+e^{-t}\prn{u(0)-u^*}$.  This linear path obviously converges to the fixed point $u^*$, and if we are able to show that the actual trajectory $u(t)$ asymptotically approaches this idealized linear path, we will have established that the system converges to $u^*$.  To take this approach, we examine the quantity
$$h(t) = u(t) - u^*-e^{-t}\prn{u(0)-u^*},$$
which is the deviation from the linear path.
Consider the norm of this deviation:
\begin{align*}
\eltwo{h(t)} &= \eltwo{u(t) - u^*-e^{-t}\prn{u(0)-u^*}} \\
 & = \eltwo{e^{-t}\dint{0}{t}{e^s\prn{\Phi^T\Phi-I} \widetilde{a}(s) ds}} \\
 & \leq e^{-t}\eltwo{\Phi^T\Phi-I}\dint{0}{t}{e^s \eltwo{\widetilde{a}(s)} ds}.
\end{align*}
To show convergence to zero, we split the integral into two parts. Since $\widetilde{a}(t)\underset{t\rightarrow+\infty}{\longrightarrow}0$, then for any $\widetilde{\epsilon}>0$, there exists a time $t_c\geq0$ such that $\forall t\geq t_c$, $\eltwo{\widetilde{a}(t)}\leq\widetilde{\epsilon}$.
Moreover, since $\eltwo{\widetilde{a}(t)}$ is continuous and goes to zero as $t$ goes to infinity, it admits a maximum $\mu$, $\forall t\in\field{R}$. This yields, for all $t\geq 2t_c$:
\begin{align*}
\eltwo{h(t)}  & \leq e^{-t}\eltwo{\Phi^T\Phi-I}\mu\dint{0}{t_c}{e^s ds} \\
 & \qquad \qquad + e^{-t}\eltwo{\Phi^T\Phi-I}\widetilde{\epsilon}\dint{t_c}{t}{e^s ds} \\
 & \leq \eltwo{\Phi^T\Phi-I} \mu \Prn{e^{-t+t_c}-e^{-t}} \\
 & \qquad \qquad + \eltwo{\Phi^T\Phi-I} \widetilde{\epsilon} \Prn{1-e^{-t+t_c}} \\
 & \leq \eltwo{\Phi^T\Phi-I} \mu \Prn{e^{-t/2}-e^{-t}} + \eltwo{\Phi^T\Phi-I} \widetilde{\epsilon} . 
\end{align*}
Since the left term converges to $0$ and $\widetilde{\epsilon}$ can be chosen arbitrarily small, this shows that the trajectory $u(t)$ converges to the trajectory $u^*+e^{-t}\prn{u(0)-u^*}$ as $t$ goes to infinity, and thus, we can conclude that
$u(t)\underset{t\rightarrow+\infty}{\longrightarrow}u^*.$

Because we have shown separately that both the active and inactive nodes converge for any initial state, it concludes our proof that the system is globally convergent.
\end{IEEEproof}

\section{Proof of Lemma \ref{lem:pties}}
\label{app:pties}

\begin{IEEEproof}
Each of the four cases will be treated separately.
\begin{enumerate}[(i)]
\item For any $\widetilde{u}_n\in\field{R}$, let $z_n=\sign{\widetilde{u}_n}$. Since the thresholding function is non-decreasing \eqref{eq:cond2}, and $T_{\lambda}(-u_n)=- T_{\lambda}(u_n)$ from \eqref{eq:cond1}, we have:
\begin{align*}
z_n=\sign{\widetilde{u}_n} & \Rightarrow 0 \leq z_n\widetilde{u}_n \\
 & \Rightarrow z_n u_n^* \leq z_n\widetilde{u}_n + z_n u_n^* \\
 & \Rightarrow T_{\lambda}\prn{z_n u_n^*} \leq T_{\lambda}\prn{z_n\widetilde{u}_n + z_n u_n^*} \\
 & \Rightarrow z_n T_{\lambda}\prn{u_n^*} \leq z_n T_{\lambda}\prn{\widetilde{u}_n + u_n^*} \\
 & \Rightarrow 0 \leq z_n \Prn{ T_{\lambda}\prn{\widetilde{u}_n + u_n^*} - T_{\lambda}\prn{u_n^*} } \\
 & \Rightarrow 0 \leq z_n \widetilde{a}_n \\
 & \Rightarrow \sign{\widetilde{a}_n} = z_n = \sign{\widetilde{u}_n}.
\end{align*}
\item We can separate this proof into four cases.

If $\abs{\widetilde{u}_n+u_n^*}\leq\lambda$ and $\abs{u_n^*}\leq\lambda$, we have:
 $\abs{\widetilde{a}_n} = \abs{T_{\lambda}(\widetilde{u}_n+u_n^*) - T_{\lambda}(u_n^*)} = 0 \leq \alpha\abs{\widetilde{u}_n}.$
If $\abs{\widetilde{u}_n+u_n^*}\leq\lambda$ and $\abs{u_n^*}>\lambda$, according to the mean value theorem, since the function $f(\cdot)$ is continuous on $\Prn{\sign{u_n^*}\lambda,u_n^*}$, differentiable on $\prn{\sign{u_n^*}\lambda,u_n^*}$ and $f(\sign{u_n^*}\lambda)=0$, there exist $c\in\prn{\sign{u_n^*}\lambda,u_n^*}$ such that:
\begin{align*}
 \abs{\widetilde{a}_n} & = \abs{f(u_n^*)} = \abs{f'(c)\prn{u_n^*-\sign{u_n^*}\lambda}} \\
 & = \abs{f'(c)}\prn{\abs{u_n^*}-\lambda} \ \ (\mathrm{since\ } \abs{u_n^*}>\lambda) \\
 & \leq \alpha \prn{\abs{u_n^*}-\lambda} \\
 & \leq \alpha \abs{\widetilde{u}_n} \ \ (\mathrm{since\ } \abs{u_n^*}-\abs{\widetilde{u}_n}\leq\abs{\widetilde{u}_n+u_n^*}\leq\lambda).
\end{align*}
If $\abs{\widetilde{u}_n+u_n^*}>\lambda$ and $\abs{u_n^*}\leq\lambda$, according to the mean value theorem, there exists $c\in\prn{\sign{\widetilde{u}_n+u_n^*}\lambda,\widetilde{u}_n+u_n^*}$ such that:
\begin{align*}
 \abs{\widetilde{a}_n} & = \abs{f(\widetilde{u}_n+u_n^*)} \\
 & = \abs{f'(c)\prn{\widetilde{u}_n+u_n^*-\sign{\widetilde{u}_n+u_n^*}\lambda}} \\
 & = \abs{f'(c)}\prn{\abs{\widetilde{u}_n+u_n^*} -\lambda} \\
 & \leq \alpha \prn{\abs{\widetilde{u}_n} + \abs{u_n^*} - \lambda} \leq \alpha \abs{\widetilde{u}_n}.
\end{align*}
Finally, if $\abs{\widetilde{u}_n+u_n^*}>\lambda$ and $\abs{u_n^*}>\lambda$, according to the mean value theorem, there exists $c\in\prn{u_n^*,\widetilde{u}_n+u_n^*}$ such that
 $\abs{\widetilde{a}_n}  = \abs{f(\widetilde{u}_n+u_n^*) - f(u_n^*)} = \abs{f'(c)\widetilde{u}_n} \leq \alpha \abs{\widetilde{u}_n}$.
\item Using property (i) and (ii), we have:
\begin{align*}
 \widetilde{a}_{\mathcal{T}}^T\widetilde{a}_{\mathcal{T}} & = \sum_{n\in{\mathcal{T}}}{\widetilde{a}_n\widetilde{a}_n} = \sum_{n\in{\mathcal{T}}}{\abs{\widetilde{a}_n}\abs{\widetilde{a}_n}} \\
 & \leq \sum_{n\in{\mathcal{T}}}{\alpha\abs{\widetilde{u}_n}\abs{\widetilde{a}_n}} = \alpha \sum_{n\in{\mathcal{T}}}{\widetilde{u}_n\widetilde{a}_n} = \alpha \widetilde{u}_{\mathcal{T}}^T\widetilde{a}_{\mathcal{T}} \\
 & \leq \sum_{n\in{\mathcal{T}}}{\alpha\abs{\widetilde{u}_n}\alpha\abs{\widetilde{u}_n}} = \alpha^2 \sum_{n\in{\mathcal{T}}}{\widetilde{u}_n\widetilde{u}_n} = \alpha^2 \widetilde{u}_{\mathcal{T}}^T\widetilde{u}_{\mathcal{T}}.
\end{align*}
\item Since the soft-threshold function is non-decreasing, then the function $g_n(s)=T_{\lambda}(s+u_n^*)-T_{\lambda}(u_n^*)$ is also non-decreasing. Moreover $g_n(0)=0$ and $g_n(\widetilde{u}_n(t))=\widetilde{a}_n(t)$. As a consequence, we can bound the integral by:
$$ \dint{0}{\widetilde{u}_n(t)}{g_n(s) ds} \leq \Bigl(g_n(\widetilde{u}_n(t))-g_n(0)\Bigr)\widetilde{u}_n(t) = \widetilde{a}_n\widetilde{u}_n .$$
As a consequence:
$$ \dsum{n=1}{N}{\dint{0}{\widetilde{u}_n(t)}{g_n(s) ds}} \leq \dsum{n=1}{N}{\widetilde{a}_n\widetilde{u}_n}  = \widetilde{a}^T\widetilde{u} .$$
\end{enumerate}
\end{IEEEproof}


\bibliographystyle{IEEEbib}
\bibliography{LcaNN}

\end{document}